\documentclass[preprint,12pt]{article}

\usepackage[english]{babel}
\usepackage{amsmath}
\usepackage{amssymb}
\usepackage{amsthm}
\usepackage{pstricks}
\usepackage{pst-plot}
\usepackage{graphicx}
\usepackage{latexsym}
\usepackage{comment}

\numberwithin{equation}{section}

\newtheorem{lemma}{Lemma}[section]
\newtheorem{theorem}[lemma]{Theorem}
\newtheorem{definition}[lemma]{Definition}

 \DeclareMathOperator{\Ai}{Ai}
 \DeclareMathOperator{\supp}{supp}
 \DeclareMathOperator{\Impart}{Im}
 \DeclareMathOperator{\Repart}{Re}
 \renewcommand{\Im}{\Impart}
 \renewcommand{\Re}{\Repart}

\title{
Asymptotic zero distribution of complex orthogonal polynomials associated with Gaussian quadrature}

\author{A. Dea\~{n}o\footnotemark[1]\ , D. Huybrechs\footnotemark[2]\ , and A.B.J. Kuijlaars\footnotemark[3]}

\date{ \ } 

\begin{document}

\maketitle
\renewcommand{\thefootnote}{\fnsymbol{footnote}}
\footnotetext[1]{Department of Mathematics, Universidad Carlos III de Madrid, Spain. \\
Email: alfredo.deanho\symbol{'100}uc3m.es.}
\footnotetext[2]{Department of Computer Science, Katholieke Universiteit Leuven, Belgium. \\
Email: daan.huybrechs\symbol{'100}cs.kuleuven.be.}
\footnotetext[3]{Department of Mathematics, Katholieke Universiteit Leuven, Belgium.  \\
Email: arno.kuijlaars\symbol{'100}wis.kuleuven.be.}

\begin{abstract}
In this paper we study the asymptotic behavior of a family of
polynomials which are orthogonal with respect to an exponential
weight on certain contours of the complex plane. The zeros of
these polynomials are the nodes for complex Gaussian quadrature of
an oscillatory integral on the real axis with a high order
stationary point, and their limit distribution is also analyzed.
We show that the zeros accumulate along a contour in the complex plane
that has the $S$-property in an external field. In addition, the
strong asymptotics of the orthogonal polynomials is  obtained by
applying the nonlinear Deift--Zhou steepest descent method to the
corresponding Riemann--Hilbert problem.
\end{abstract}


\section{Introduction}

\subsection{Oscillatory integrals}

We study the limiting behavior of the zeros of the polynomials
that are orthogonal with respect to an oscillatory weight function
of exponential type along a path $\Gamma$ in the complex plane,
\begin{equation}\label{E:orthogonality}
    \int_{\Gamma} \pi_n(z)z^k e^{i z^r}d{z}=0, \qquad k=0,1,\ldots,n-1,
\end{equation}
such that the integral is well defined. Parameter $r \geq 2$ is an
integer and we will focus mainly on the case $r=3$.

Our motivation originates in a Fourier-type integral on a finite interval of the real axis of the general form
\begin{equation}\label{E:I_general}
I[f] = \int_a^b f(x) e^{i \omega g(x)} d{x},
\end{equation}
where $\omega > 0$ is a frequency parameter, $f$ is called the \emph{amplitude} and $g$ is the \emph{phase} or \emph{oscillator}. Integrals of this type appear in many scientific disciplines involving wave phenomena, such as acoustics, electromagnetics and optics (see for example \cite{huybrechs:2009:hoq} and references therein). For $\omega\gg 1$, integrals of this kind are a recurring topic in asymptotic analysis, and we recall in particular the classical method of steepest descent, which can be applied when $f$ and $g$ are analytic in a neighbourhood of $[a,b]$, see for instance \cite{wong:2001:asymptotic}.

We will concentrate on the case where the oscillator $g$ has a single stationary point $\xi$ of order $r-1$ inside the interval $[a,b]$, with $r\geq 2$, i.e., $g^{(j)}(\xi)=0$, $j=1,\ldots,r-1$ but $g^{(r)}(\xi) \neq 0$. Without loss of generality, we take this point $\xi$ to be the origin and the canonical example is the following:
\begin{equation}\label{E:I_canonical}
 I[f] := \int_a^b f(x) e^{i \omega x^{r}} d{x},
\end{equation}
with $a < 0$ and $b > 0$. Assuming a single stationary point, the general form~\eqref{E:I_general} can always be brought into this form by a change of variables.

\subsection{Numerical evaluation}

\begin{figure}[t]
\begin{center}
\includegraphics{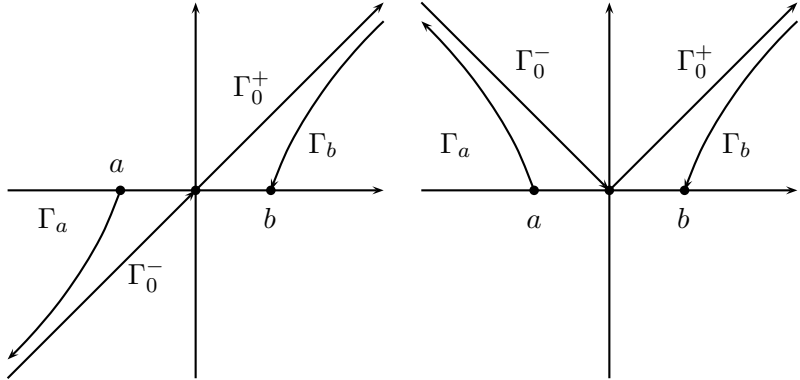}
\caption{Approximate contours of integration in the complex plane
corresponding to even $r$ (left) and odd $r$
(right).}\label{F:saddlepoint_contours}
\end{center}
\end{figure}

We assume $f$ analytic in a complex neighbourhood of the interval
$[a,b]$. As shown in \cite{Huybrechs:06:osc1}, one possible
numerical strategy for the evaluation of~\eqref{E:I_canonical} is
to consider paths of steepest descent stemming from the endpoints
and from the stationary point. In this way, we can decompose the
original integral as follows:
$$
\int_a^b f(x)e^{i\omega x}d x=
\left(\int_{\Gamma_a}+\int_{\Gamma^{-}_0}+\int_{\Gamma^{+}_0}+\int_{\Gamma_b}\right)f(x)e^{i\omega x}d x,
$$
where the paths are depicted in Fig.~\ref{F:saddlepoint_contours}.
Making an appropriate change of variables, the line integrals
along these paths have the form
\begin{equation*}
 \int_0^\infty u(z) e^{-\omega z^\mu} d{z},
\end{equation*}
with $\mu=1$ for the endpoints and $\mu=r$ for a stationary
point of order $r-1$. Each of these integrals can be efficiently approximated
using Gaussian quadrature, because the optimal
polynomial order of Gaussian quadrature translates into optimal
asymptotic order in this setting: for $n$ quadrature points the
error behaves like $\mathcal{O}(\omega^{-\frac{2n+1}{\mu}})$ as
$\omega \to \infty$, see \cite{DH:2008:CG}. This order is
approximately twice that of a classical asymptotic expansion
truncated after $n$ terms.

There are two paths of steepest descent originating from the
stationary point, called $\Gamma^{-}_0$ and $\Gamma^{+}_0$ in
Fig.~\ref{F:saddlepoint_contours}. Both paths are straight lines
and their structure depends essentially on the parity of $r$: for
odd $r$, these lines form an angle equal to $\pi-\pi/r$, whereas
for even $r$ they form one straight line in the complex plane.

In order to keep the total number of function evaluations to a
minimum, it is desirable to evaluate both line integrals using only one
quadrature rule of Gaussian type~\cite{DH:2008:CG}. This amounts
to constructing a quadrature rule for the functional
\begin{equation} \label{E:Mf}
 M[f] = \int_{\Gamma} f(z) e^{i z^r} d{z},
\end{equation}
where $\Gamma=\Gamma^{-}_0\cup\Gamma^{+}_0$ is the concatenation of the two steepest descent paths through  the origin.

In the case $r=2$, this leads to classical Gauss-Hermite
quadrature, which involves the weight function $e^{-x^2}$ on the
real line $(-\infty,\infty)$. Higher even values of $r$ lead to
straightforward generalizations, and in all cases the quadrature points lie on the paths of steepest descent.

\begin{figure}[t]
\centerline{\includegraphics[height=65mm,width=160mm]{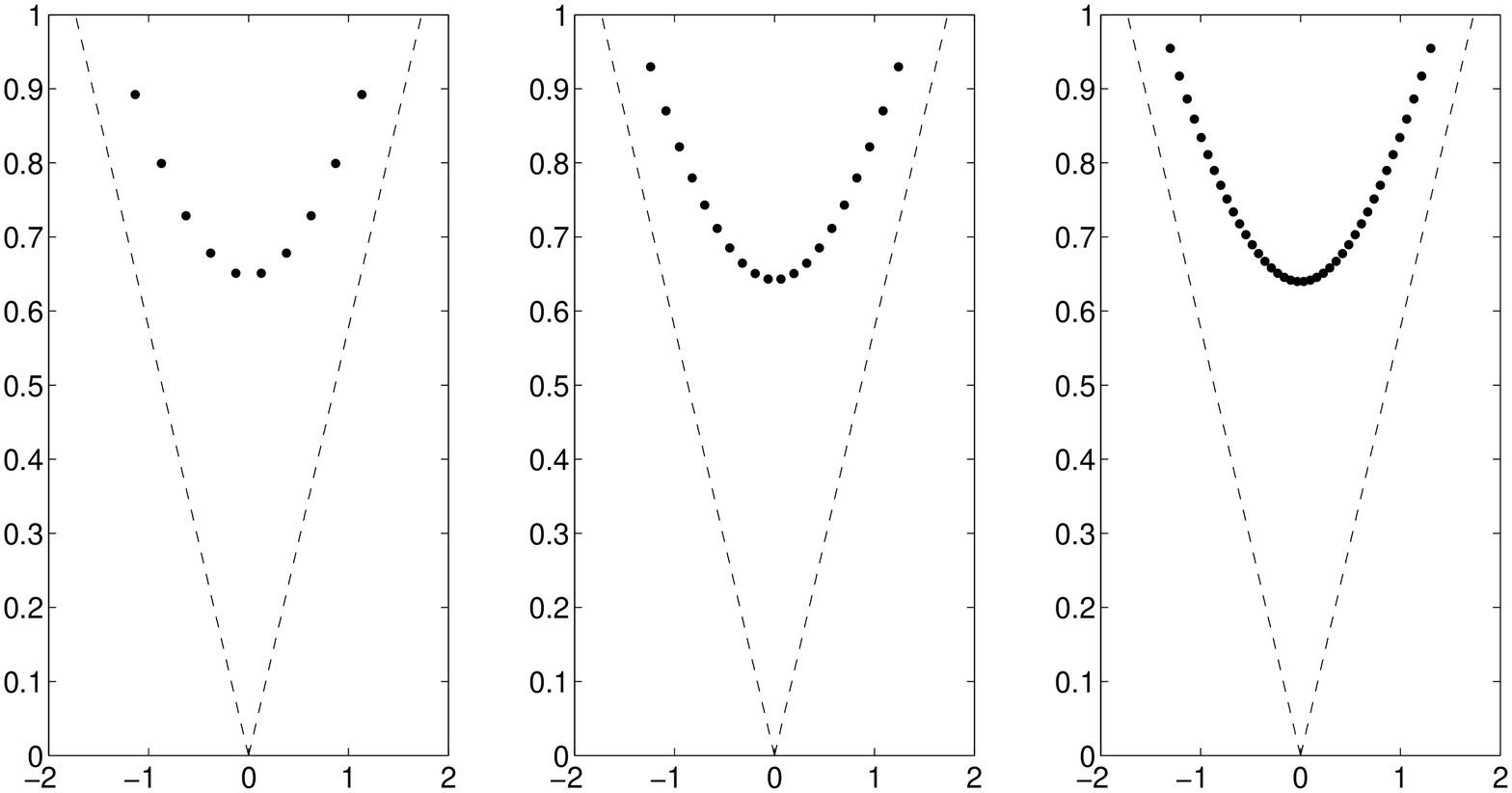}}
\caption{Location of the quadrature nodes for $r=3$ on $[-1,1]$,
corresponding to $n=10$ (left), $n=20$ (center) and $n=40$
(right). In dashed line, the paths of steepest descent from the
origin.}\label{F:indefinite_r3}
\end{figure}

For odd $r$, the functional \eqref{E:Mf} is indefinite and the
existence of orthogonal polynomials is not guaranteed a priori.
Nevertheless, the orthogonal polynomials and their zeros can be
computed numerically. However, one finds that the zeros, which are
the complex quadrature points for the integral (3), do not lie on
the paths of steepest descent anymore. Instead, they seem to lie
on a curve in a sector of the complex plane bounded by the paths
of steepest descent. Their location for $r=3$ is shown in
Fig.~\ref{F:indefinite_r3} for several values of $n$. Similar
phenomena are observed for larger odd values of $r$.

\subsection{Orthogonality in the complex plane}

The problem of Gaussian quadrature leads to the study of
polynomials $\pi_n(z)$ that are orthogonal in the sense of
\eqref{E:orthogonality}, where $r$ is a positive integer ($r\geq
3$ for a non-classical case) and $\Gamma$ is the combination of two
paths of steepest descent of the exponential function $e^{i
z^r}$ from the origin, so
$$
\arg z=\frac{\pi}{2r}+\frac{2m\pi}{r}, \qquad m=0,1,\ldots,r-1.
$$
The straight lines $\Gamma^{-}_0$ and $\Gamma^{+}_0$ in
Fig.~\ref{F:indefinite_r3} correspond respectively to $m=0$ and
$m=\lfloor \frac{r}{2} \rfloor$. In the case $r=3$, these lines
form angles of $\pi/6$ and $5\pi/6$ with respect to the positive
real axis.

Putting $\lambda_n = (n/r)^{1/r}$ and
\begin{equation} \label{E:Pn}
    P_n(z) = \lambda_n^{-n} \pi_n(\lambda_n z)
\end{equation}
we note that~\eqref{E:orthogonality} can be written in the form
\begin{equation}
\label{E:varying_weight}
  \int_{\Gamma} P_n(z) z^k e^{-nV(z)}d z=0, \qquad k=0, \ldots, n-1,
\end{equation}
where
\begin{equation}\label{E:V}
V(z)=-i z^r/r.
\end{equation}
Note that \eqref{E:Pn} is again a monic polynomial, and the zeros of
$P_n(z)$ and $\pi_n(z)$ are the same but for rescaling with the
parameter $\lambda_n$.

The orthogonality \eqref{E:varying_weight} is an example of
non-Hermitian orthogonality with respect to a varying weight on a
curve in the complex plane. A basic observation is that the path
$\Gamma$ of the integral in \eqref{E:varying_weight} can be
deformed into any other curve that is homotopic to it in the
finite plane, and that connects the same two sectors at infinity. For
any such deformed $\Gamma$ we still have the orthogonality
condition \eqref{E:varying_weight}.


In order to find where the zeros of $P_n(z)$ lie for large $n$, we
have to select the `right' contour.  Stahl \cite{Stahl:1986} and
Gonchar--Rakhmanov \cite[Sec.~3]{GR:1989:eq} studied and solved this problem, and from their works it is known that the appropriate contour should have a
symmetry property (the so-called $S$-property) in the sense of logarithmic
potential theory with external fields.
We recall this concept in the next subsection.

In the case \eqref{E:V} with $r=3$, we can identify the
curve with the $S$-property explicitly as a critical trajectory
of a quadratic differential. Other cases where the potential problem is explicitly solved include
\cite{GR:1989:eq}, \cite{MF:1997} and \cite{Apt:2002} in connection
with best rational approximation of $e^{-x}$ on $[0,\infty)$,
\cite{Baik:2001:random}  in connection with a last passage percolation problem, and \cite{KML:2001:Laguerre}, \cite{KML:2004:Laguerre}, \cite{KMF:2004:Jacobi}, \cite{MFMGO:2001:Laguerre}, \cite{MFO:2005:Jacobi} in connection with  classical orthogonal polynomials (Laguerre and Jacobi) with non-standard parameters.

Varying orthogonality on complex curves is treated in
detail in the more recent accounts \cite{Apt:2007:complex}
and \cite{BertolaMo:2009}, which contain a Riemann-Hilbert
steepest descent analysis in a fairly general setting, assuming
the knowledge of the curve with the $S$-property. See also
\cite{Bertola:Boutroux} for an approach based on algebraic
geometry and Boutroux curves and \cite{Bertola:Lp} for
extensions to $L_p$ optimal polynomials.

\subsection{The $S$ property}


Let $V$ be a polynomial.
We consider a smooth curve $\Gamma\subset\mathbb{C}$, such that the integral
in \eqref{E:varying_weight} is well-defined, and we want to minimize the
weighted energy:
\begin{equation} \label{E:energyonGamma}
    I_V(\nu)
    = \iint \log\frac{1}{|z-s|}d\nu(z)d\nu(s)+ \Re \int
    V(s)d\nu(s),
\end{equation}
among all Borel probability measures $\nu$ supported on $\Gamma$.
Following the general theory of logarithmic potential theory with
external fields, see \cite{ST:1997:Pot}, this problem
has a unique solution, which is called the equilibrium measure on
$\Gamma$ in the presence of the external field $\Re V$. We denote this equilibrium measure by $\mu$.

Let
 \begin{equation} \label{E:Umupotential}
   U^{\mu}(z) = \int \log \frac{1}{|z-s|} d\mu(s)
\end{equation}
 be the logarithmic potential of $\mu$. It satisfies
\begin{equation} \label{E:Uequilibrium}
\begin{aligned}
  2 U^{\mu}(z) + \Re V(z) & = \ell, \qquad  z\in \supp \mu, \\
  2 U^{\mu}(z) + \Re V(z) & \geq \ell, \qquad z \in \Gamma \setminus \supp \mu,
\end{aligned}
\end{equation}
for some constant $\ell$, see  \cite{ST:1997:Pot}. If $\Gamma$ is an analytic contour, then $\supp \mu$ will consist
of a finite union of analytic arcs. Now we can define the $S$-property.

\begin{definition}
The analytic contour $\Gamma$ has the $S$-property in the external field
$\Re V$ if for every $z$ in the interior of the analytic arcs that
constitute $\supp \mu$, we have
\begin{align} \label{E:Sproperty}
   \frac{\partial}{\partial n_+}  \left[ 2 U^{\mu}(z) + \Re V(z) \right] =
    \frac{\partial}{\partial n_-} \left[ 2 U^{\mu}(z) + \Re V(z) \right].
\end{align}
Here $\frac{\partial}{\partial n_{\pm}}$ denote the two normal
derivatives taken on either side of $\Gamma$.
\end{definition}

The result of Gonchar-Rakhmanov then reads (for the special
case of polynomial $V$):

\begin{theorem} \textbf{Gonchar-Rakhmanov \cite[Sec.~3]{GR:1989:eq}} \label{th:Sproperty}
If $\Gamma$ is a contour with the $S$-property \eqref{E:Sproperty}
in the external field $\Re V$, then the equilibrium measure $\mu$ on $\Gamma$ in the
external field $\Re V$ is the weak limit of the
normalized zero counting measures of the polynomials $P_n$
defined by the orthogonality \eqref{E:varying_weight}.
\end{theorem}

\subsection{Outline of the paper}
In the next section we present the main results of this paper,
corresponding to \eqref{E:V} with $r=3$. These can be summarized in the following points:
\begin{itemize}
\item We present a finite curve $\gamma\subset\mathbb{C}$, which
is a critical trajectory of a certain quadratic differential $Q(z)d z^2$,
see Theorem \ref{theoremgamma}.
\item We prove that this curve $\gamma$ can be prolonged to $\infty$ in a
suitable way, thus obtaining a curve $\Gamma$ with the $S$-property in the
presence of the external field $\Re V$, see Theorem \ref{theo:Sproperty}.
\item As a consequence of Gonchar-Rakhmanov theorem, it is possible to obtain the weak limit distribution of the zeros of $P_n(z)$  as $n\to\infty$, see Theorem \ref{theoremr3}.
\item Additionally, a full Riemann--Hilbert analysis of this problem is feasible and yields both existence of the sequence of orthogonal polynomials $P_n(z)$ for large enough $n$ and the asymptotic behavior of $P_n(z)$ in various regions of the complex plane as $n\to\infty$, see Theorem \ref{th:strongasymptotics}.
\end{itemize}

\section{Statement of results}

\subsection{Definition of the curve $\gamma$}
In the case \eqref{E:V} with $r=3$, the curve with the $S$-property  is given in terms of the critical trajectory of the
quadratic differential $Q(z)d z^2$, where
\begin{equation} \label{Qr3}
Q(z)=-\frac 14(z+i)^2 (z^2-2i z-3).
\end{equation}

The polynomial \eqref{Qr3} has a double root at $z=-i$ and two
simple roots at $z_1 = -\sqrt{2} + i$ and $z_2 = \sqrt{2} +
i$.

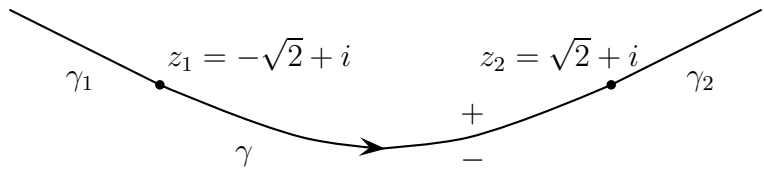
\begin{figure}[t]
\begin{pspicture}(0,0)(10,3.5)
\pscurve{->,arrowsize=0.25}(3,1.5)(4.9,0.8)(6,0.65)
\pscurve(5.9,0.65)(7.1,0.8)(9,1.5)
\psline(1,2.5)(3,1.5)
\psline(9,1.5)(11,2.5)
\psdot(3,1.5)
\psdot(9,1.5)
\put(3.1,1.75){$z_1 = -\sqrt{2} + i$}
\put(7.25,1.75){$z_2 = \sqrt{2} + i$}
\put(4,0.5){$\gamma$}
\put(7,1){$+$}
\put(7,0.4){$-$} \put(1.75,1.5){$\gamma_1$}
\put(10,1.5){$\gamma_2$}
\end{pspicture}
\caption{The contour $\Gamma$ consists of the
critical trajectory $\gamma$ and its analytic
continuations $\gamma_1$ and $\gamma_2$.  \label{Figure3}}
\end{figure}

The critical trajectory $\gamma$ is an analytic arc from
$z_1$ to $z_2$ so that
\begin{equation} \label{E:Dz}
    \frac{1}{\pi i} \int_{z_1}^z Q^{1/2}(s) d s
\end{equation}
is real for $z\in\gamma$, see \cite{strebel:1984:quad}. We first show that this curve indeed exists.

\begin{theorem} \label{theoremgamma}
There exists a critical trajectory $\gamma$ of
the quadratic differential $Q(z)d z^2$, where $Q(z)$ is given in \eqref{Qr3},
that connects the two zeros $z_1=-\sqrt{2}+i$ and $z_2=\sqrt{2}+i$ of $Q$.
\end{theorem}

The proof of the theorem is contained in Section \ref{s:proof1}.

\subsection{Contour with $S$-property}
In what follows we use the analytic arc $\gamma$, whose existence is guaranteed
by Theorem \ref{theoremgamma}, with an orientation so that $z_1$ is the starting point of $\gamma$
and $z_2$ is the ending point. The $+$ side ($-$ side) of $\gamma$ is on the left (right)
as we traverse $\gamma$ according to its orientation, as shown in Figure~\ref{Figure3}.

From now on the square root $Q^{1/2}(z)$ is defined with
a branch cut along $\gamma$ and so that
    \begin{equation} \label{E:Qsqrt}
    Q^{1/2}(z) = - \frac{1}{2} i z^2  - \frac{1}{z} + \mathcal{O}\left(\frac{1}{z^2}\right)
    \end{equation}
as $z \to \infty$. This branch is then used for
example in \eqref{E:Dz}.  We use $Q^{1/2}_+(s)$, when $s \in \gamma$, to denote the limiting value of $Q^{1/2}(z)$ as $z$ approaches $s \in \gamma$
from the $+$ side.

The curve $\gamma$ has an analytic extension to an
unbounded oriented contour
\begin{equation} \label{E:Gamma}
    \Gamma = \gamma_1\cup\gamma\cup\gamma_2
    \end{equation}
that we  use for the orthogonality \eqref{E:varying_weight}.
The parts $\gamma_1$ and $\gamma_2$ are such that
\begin{equation} \label{E:isreal1}
   \phi_1(z) = \int_{z_1}^z Q^{1/2}(s) ds \quad \text{ is real and positive  for } z \in \gamma_1,
    \end{equation}
and
\begin{equation} \label{E:isreal2}
    \phi_2(z) = \int_{z_2}^z Q^{1/2}(s) ds \quad \text{ is real and positive for } z \in \gamma_2.
    \end{equation}

The main result of the paper is then the following.

\begin{theorem} \label{theo:Sproperty}
The contour $\Gamma$ is a curve with the $S$-property in
the external field $\Re V$. In addition, we have that the equilibrium measure on $\Gamma$ in the external field $\Re V$ is given by
the probability measure
\begin{equation} \label{E:mu}
    d\mu(s) = \frac{1}{\pi i}\, Q^{1/2}_+(s) d s.
\end{equation}
\end{theorem}


The proof of this theorem is presented in Section \ref{s:proof2}.
The general result of Gonchar-Rakhmanov, see Theorem \ref{th:Sproperty}, then implies:

\begin{theorem} \label{theoremr3}
Assume $r=3$. For large enough $n$ the monic polynomial $P_n(z)$ of degree $n$
satisfying the orthogonality relation \eqref{E:varying_weight}
with $V(z)$ given by \eqref{E:V} exists uniquely. Furthermore, denoting by
$$
   z_{1,n}, \ldots, z_{n,n}
$$
its $n$ zeros in the complex plane, we have
\begin{enumerate}
\item[\rm (a)] as $n\to\infty$ the zeros accumulate on $\gamma$;
\item[\rm (b)] the normalized zero
counting measures have a weak limit
\[ \frac 1n \sum_{k=1}^n \delta_{z_{k,n}}
\stackrel{*}{\longrightarrow} d \mu,
\]
where $\mu$ is given by \eqref{E:mu}.
\end{enumerate}
\end{theorem}

\subsection{Riemann-Hilbert analysis}

Theorem \ref{theoremr3} also follows from a steepest descent analysis
for the Riemann-Hilbert problem that characterizes the polynomials $P_n$.
We include an exposition of this method in Section \ref{s:proof4} for two main reasons:

\begin{itemize}
\item The Riemann-Hilbert analysis not only provides the limit behavior
of the distribution of zeros of $P_n(z)$, but also
strong asymptotics of the orthogonal polynomials in the complex plane.
\item The present case can be viewed as a simple model problem
for Riemann-Hilbert analysis in the complex plane. We hope that the
present analysis can be useful as an introduction to this
powerful method.
\end{itemize}

The Riemann-Hilbert problem
for orthogonal polynomials was found by  Fokas, Its and Kitaev
\cite{fokas:1992:isomonodromy}, and the steepest descent analysis of Riemann-Hilbert problems is due
to Deift and Zhou \cite{deift:1993:steepestdescent}.  The steepest descent analysis
for orthogonal polynomials with varying weights on the real line
is due to Deift et al., see \cite{DKMVZ:1999:varying} and \cite{Deift:2000:RH}.

The extension of this method to orthogonal polynomials on curves
in the complex plane is not new. It has already been presented in various
papers, see for example \cite{Apt:2002},
\cite{Apt:2007:complex}, \cite{BertolaMo:2009}, \cite{DKMVZ:1999:varying} and \cite{KML:2001:Laguerre}.
However, an attractive feature of the example treated here is
that all quantities in the analysis can be computed explicitly.
In that respect it is similar to \cite{KML:2001:Laguerre}.

In order to formulate the additional asymptotic results that
follow from the steepest descent analysis, we introduce some
more notation. We use the function $\phi_2(z)$ defined in
\eqref{E:isreal2} and the related function (the $g$-function)
\begin{equation} \label{E:defg}
    g(z) =\frac 12 V(z)-\phi_2(z)-l, \qquad l = \frac{1}{3} + \frac{1}{2} \log 2,
    \end{equation}
which has an alternative expression
\[ g(z) = \int \log(z-s) d\mu(s) \]
in terms of the equilibrium measure $\mu$ on $\gamma$.
In the case $r=3$ there is an explicit expression for $\phi_2(z)$:
\begin{multline} \label{E:phi_explicit}
    \phi_2(z) = -\frac{i}{6}z(z+i)\sqrt{z^2-2i z-3} \\ -
        \log(z-i+ \sqrt{z^2-2i z-3})+\frac{1}{2} \log 2.
\end{multline}

The so-called global parametrix $N(z)$ is defined in terms of the function
\begin{equation} \label{E:beta}
 \beta(z) = \left(\frac{z-z_2}{z-z_1}\right)^{1/4}, \qquad z \in \mathbb C \setminus \gamma,
\end{equation}
with the branch cut taken along $\gamma$. $N(z)$ is a $2 \times 2$ matrix
valued function with entries
$$
 N_{11}(z) = N_{22}(z) = \frac{\beta(z)+\beta(z)^{-1}}{2}, \qquad
 N_{12}(z)= - N_{21}(z) = \frac{\beta(z)-\beta(z)^{-1}}{2i},
$$
that also appear in the asymptotic formulas in Theorem \ref{th:strongasymptotics}.

Finally, near the endpoint $z_2$ we require a conformal map
\begin{equation}
      f(z) = \left[ \frac{3}{2} \phi_2(z) \right]^{2/3},
\end{equation}
which maps $\gamma$ and $\gamma_2$ near $z_2$ into the real line.
A local Riemann-Hilbert problem is solvable
explicitly in terms of the usual Airy function $\Ai(z)$ and its
derivative $\Ai'(z)$. These functions appear in the asymptotic
formula in part (c) of Theorem \ref{th:strongasymptotics} that is
valid in a neighborhood of $z_2$.

The steepest descent analysis of the Riemann-Hilbert problem then
leads to the following result:
\begin{theorem} \label{th:strongasymptotics}
Assume $r=3$. Let $U_{\delta}(z_1)$
and $U_{\delta}(z_2)$ be small  neighbourhoods of the points $z_1$
and $z_2$ given before. As $n\to\infty$ the polynomial $P_n(z)$
has the following asymptotic behavior:
\begin{itemize}
\item[\rm (a)] Uniformly for $z$ in compact subsets of $\overline{\mathbb C}
\setminus \gamma$, we have
 \begin{equation} \label{E:asymptotic outside}
   P_n(z) = N_{11}(z) e^{n g(z)} \left(1 + O(1/n)\right)
 \end{equation}
 as $n \to \infty$;
\item[\rm (b)] There is a neighbourhood $U$ of $\gamma$ in the complex
plane, so that, uniformly for $z \in U \setminus (U_{\delta}(z_1)\cup U_{\delta}(z_2))$:
\begin{equation} \label{E:asymptotic_away}
P_n(z) = e^{n\left[\frac {V(z)}{2}-l\right]}
\left(e^{-n\phi_2(z)}N_{11}(z) \pm i e^{n\phi_2(z)}N_{12}(z) +
O(1/n) \right),
\end{equation}
where the $+$ ($-$) sign in \eqref{E:asymptotic_away} is valid for $z$ in the part of
$U$ that lies above (below) the curve $\gamma$;
\item[\rm (c)] Uniformly for $z\in U_{\delta}(z_2)$, we have as $n \to \infty$,
\begin{align*}
P_n(z) & = \sqrt{\pi}e^{n\left[\frac{V(z)}{2} - l \right]}
    \left( n^{1/6} f^{1/4}(z) \beta^{-1}(z)\Ai(n^{2/3}f(z))\left(1 + \mathcal{O}(1/n)\right) \right. \\
      & \left. \hspace*{25mm} - n^{-1/6} f^{-1/4}(z) \beta(z) \Ai'(n^{2/3} f(z)) \left(1 + \mathcal{O}(1/n)\right)\right),
\end{align*}
with the same constant $l$ as in \eqref{E:defg}.
\end{itemize}
\end{theorem}

\section{Proof of Theorem \ref{theoremgamma}}\label{s:proof1}

It follows from the general theory, see \cite{strebel:1984:quad}, that
three trajectories of the quadratic differential $Q(z)d z^2$ emanate from
each simple zero of $Q$. The three trajectories through $z_1$ emanate from
$z_1$ at angles $\theta$ that satisfy
\[ 3 \theta = \pi - Q'(z_1)  \quad (\textrm{mod}\, 2\pi). \]
From the explicit formula of $Q$ and $z_1$ we find $Q'(z_1) = - \sqrt{2} - 4 i$
and the three angles at $z_1$ are
\[ \theta =  -\frac{1}{3} \arctan(2 \sqrt{2}) + \frac{2 k\pi}{3}, \qquad k=0,1,2. \]
We let $\gamma$ be the trajectory that emanates from $z_1$ at angle
\[ \theta_0 = -\frac{1}{3} \arctan(2\sqrt{2}) = -0.4103\cdots. \]

\begin{figure}[t]
\centerline{\includegraphics[height=65mm,width=90mm]{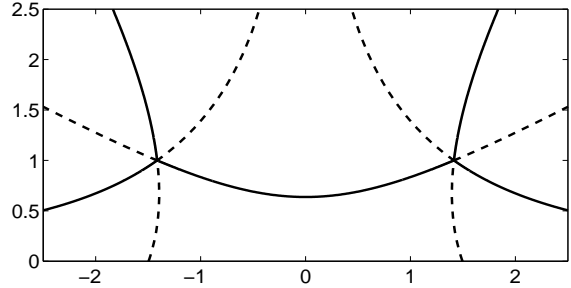}}
\vspace{-15mm}
\caption{Contour lines of $D(z)$ in the cubic case.
In solid line $\Im D(z)=0$, and in dashed line, $\Re D(z)=0$ and
$\Re D(z)= 0$ (left and right of the figure respectively).}
\label{levelcurves_r3}
\end{figure}

Let
\[ D(z) = \frac{1}{\pi i} \phi_1(z) = \frac{1}{\pi i} \int_{z_1}^z Q^{1/2}(s) d s,
    \qquad z \in \mathbb C \setminus \gamma. \]

Figure \ref{levelcurves_r3} shows the level curves $\Re
D(z)=0$, $\Re D(z)=1$, and $\Im D(z)=0$. In order to prove that
$z_1$ and $z_2$ are indeed connected by $\gamma$ (as suggested by the figure),
we use arclength parametrization of $\gamma$
\[ \gamma: \quad z = z(t), \qquad z(0) = z_1. \]
Then
$$ \int_{z(0)}^{z(t)} Q^{1/2}(s)d s = \pi i f(t),
$$
where $f(t)$ is real. Differentiating and squaring, we obtain
$$
Q(z(t))[z'(t)]^2=-\pi^2 [f'(t)]^2,
$$
which implies that
$$
\arg Q(z(t))+2\arg z'(t)=\pi \quad (\textrm{mod}\, 2\pi).
$$

\begin{figure}[t]
\centerline{\includegraphics[height=65mm,width=90mm]{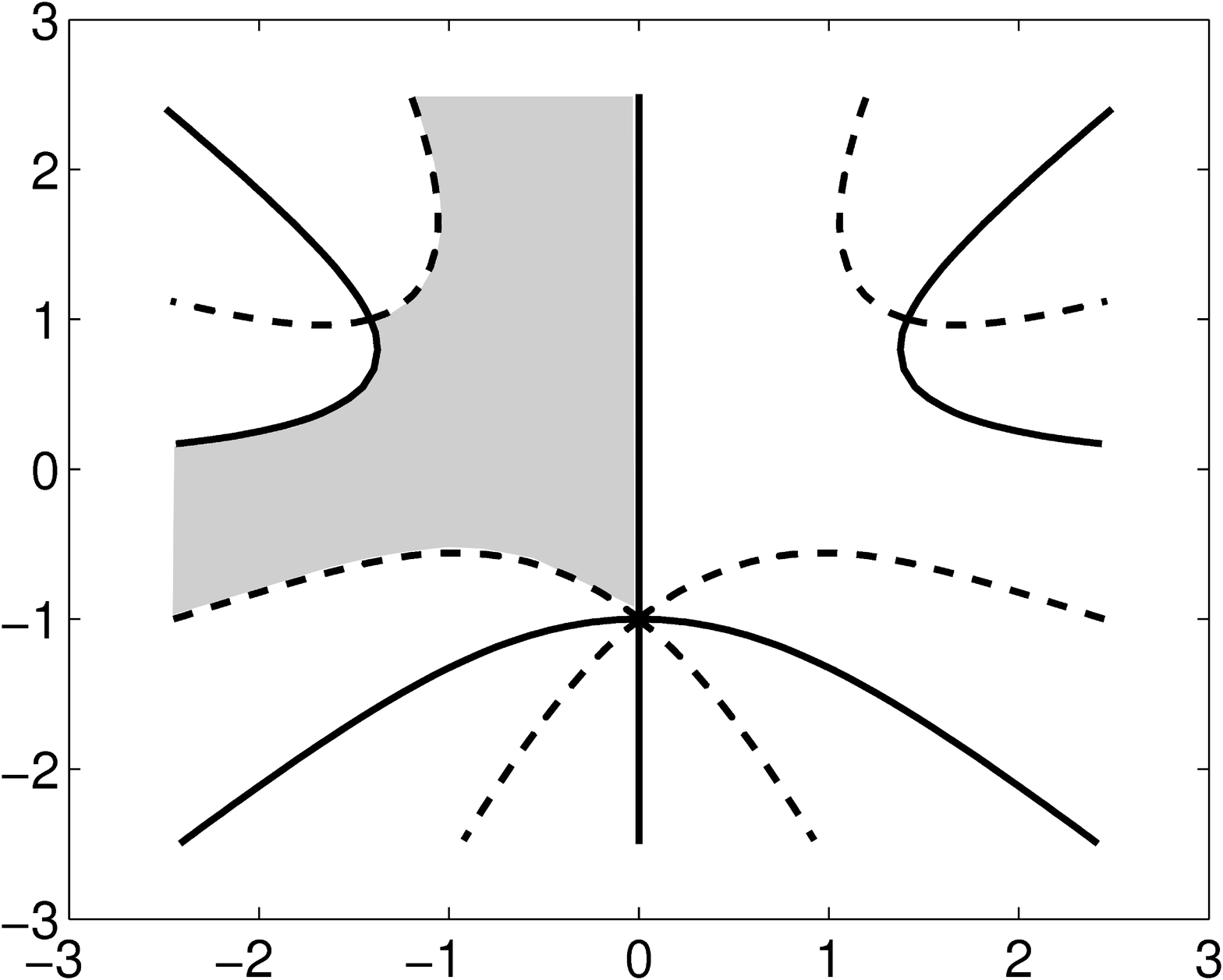}}
\caption{Contour lines of $Q(z)$ in the cubic case. In solid line
$\Re Q(z)=0$, and in dashed line $\Im Q(z)=0$.}
\label{levelcurves_Q}
\end{figure}

The level lines $\Re Q(z) = 0$ and $\Im Q(z) = 0$
are shown in Fig.~\ref{levelcurves_Q}. The two level lines
intersect of course at the zeros $z_1$, $z_2$ and $-i$ of $Q$.
The critical trajectory $\gamma$ starts at $z_1$ at an
angle $\theta_0 = -\frac{1}{3} \arctan(2\sqrt{2})$, and
therefore $\gamma$ enters the shaded region of Fig.~\ref{levelcurves_Q}, which is the region where $\Re Q(z)<0$ and $\Im Q(z)<0$,
hence $-\pi<\arg Q(z)<-\pi/2$. As a consequence
we have that
$$
-\frac{\pi}{4}<\arg z'(t)<0
$$
as long as $z(t)$ is in the shaded region. This implies that
the real part of $z(t)$ increases faster than the imaginary
part decreases. It follows that the part of $\gamma$ that
is in the shaded region is contained in the triangle with
vertices $z_1=-\sqrt{2}+i$, $i$ and $(1-\sqrt{2})i$.
Hence $\gamma$ leaves the shaded region at a point
on the imaginary axis above the other critical
point $z_0=-i$. Then by the symmetry with respect to the imaginary
axis we conclude that $\gamma$ indeed connects $z_1$ and
$z_2$.

This completes the proof of Theorem \ref{theoremgamma}.

\section{Proof of Theorem \ref{theo:Sproperty}}\label{s:proof2}
We start by presenting a characterization of the curve $\Gamma$ that is equivalent to the $S$-property.

For complex $z$, we define the $g$-function
\begin{equation} \label{E:gfunction}
    g(z)=\int_{\Gamma} \log(z-s) d \mu(s),
\end{equation}
which is analytic when $z\in\mathbb{C}\setminus\Gamma$. We observe
that $\Re g(z)=-U^{\mu}(z)$, where $U^{\mu}$ is the logarithmic
potential as defined in \eqref{E:Umupotential}. The equilibrium
properties \eqref{E:Uequilibrium} then translate into
\begin{equation} \label{E:gequilibrium}
\begin{aligned}
    \Re(-g_+(z)-g_-(z)+V(z)) & = \ell,  \qquad z\in \supp \mu,\\
    \Re(-g_+(z)-g_-(z)+V(z)) & \geq \ell,  \qquad z\in\Gamma\setminus\supp \mu.
\end{aligned}
\end{equation}

Let us write $\gamma = \supp \mu$, then from the Cauchy-Riemann equations it follows that the $S$-property
\eqref{E:Sproperty} is equivalent to the property that the
imaginary part of $-g_+ - g_- + V$ is locally constant on $\gamma$,
that is
\begin{align} \label{E:gSproperty}
    \Im (-g_+(z)-g_-(z)+V(z)) = \tilde{\ell}, \qquad z \in \gamma
\end{align}
with a possibly different constant $\tilde{\ell}$ on the different
components of $\gamma$. Then as a consequence of \eqref{E:gequilibrium}
and \eqref{E:gSproperty} we have that
\begin{equation} \label{E:constantl}
    -g_+(z)-g_-(z)+V(z) = \ell+ i\tilde{\ell}
\end{equation}
is constant on each connected component of $\gamma$.
Differentiating \eqref{E:constantl} we obtain
\begin{equation}\label{diffg}
    -g'_+(z) - g'_-(z) + V'(z) = 0, \qquad z\in \gamma.
\end{equation}

Next we observe that the function $\tfrac{1}{2} V'(z)-g'(z)$ is
analytic for $z\in \mathbb C \setminus \gamma$, and furthermore using
\eqref{diffg}:
\begin{equation*}
(\tfrac 12 V'(z)-g'(z))_+=-\tfrac 12 V'(z)+g'(z)_-
=-(\tfrac 12 V'(z)-g'(z))_-
\end{equation*}
for $z \in \gamma$.
Hence  $\tfrac 12 V'(z)-g'(z)$ has a multiplicative jump of $-1$ on $\gamma$, and therefore
\begin{equation} \label{E:Qdefinition}
    Q(z):=\left(\tfrac 12 V'(z)-g'(z)\right)^2
\end{equation}
is analytic in the whole complex plane. The asymptotic behavior
of $Q(z)$ for $z \to\infty$ follows from the fact that
$V(z)$ is a polynomial  and
\begin{equation} \label{E:gprime}
    g'(z)=\int \frac{1}{z-s} d\mu(s) = \frac{1}{z}
        + \mathcal{O} \left( \frac{1}{z^2} \right), \qquad z \to \infty,
\end{equation}
since $\mu$ is a probability measure on $\gamma$. The general case
of Liouville's theorem implies that $Q$ is a polynomial of degree
$2r-2$ if $\deg V = r$.

In general this is not enough to determine the curve $\gamma$, and we need more information on the roots of $Q(z)$ (or extra assumptions, such as that we are in the one-cut case). Since $Q^{1/2}(z)$ is analytic in $\mathbb C \setminus \gamma$
we can deduce that any zero of odd multiplicity of $Q$ is in $\gamma$. Zeros of even multiplicity can be anywhere and are typically not in $\gamma$.

Let $z_1$ be a zero of $Q$ of odd multiplicity.
From \eqref{E:gprime} we see that for $z \in \gamma$,
in the same connected component as $z_1$, we have
\[ \int_{z_1}^z Q^{1/2}_+(s) d s \in i \mathbb R. \]
This is the condition that characterizes a trajectory of
the quadratic differential $Q(z) d z^2$, emanating from a zero $z_1$ of $Q$.

In the case $V(z)=-iz^3/3$, it follows from \eqref{E:Qdefinition}
and \eqref{E:gprime} that $Q$ should be taken as a polynomial
of degree $4$ so that
\begin{equation} \label{E:QformulaC}
    Q(z) = \left(-\frac{i z^2}{2}-\frac{1}{z} +
    \mathcal{O}\left(\frac{1}{z^2}\right)\right)^2 = -\frac{z^4}{4}+i z + C,
\end{equation}
where $C$ needs to be determined. In order to do this, we make the assumption (to be justified later)
that we are in the one-cut case, that
is we assume that $\gamma$ is a single curve.  The endpoints of
the curve are then simple zeros of $Q$. Since $Q$ has degree four
there are two more zeros, which in the one-cut case, should
combine into a double zero.

In our case, there is a symmetry about the imaginary axis, and
therefore the double root should be on the imaginary axis, say $z=z_0$,
and two simple roots are symmetric with respect to the imaginary
axis, say $z_1$ and $z_2=-\overline{z}_1$. This leads to
$$
    Q(z)=-\frac 14(z-z_0)^2 (z-z_1)(z+\bar{z}_1),
$$
 which combined with \eqref{E:QformulaC} yields
\[ z_0=-i, \qquad \text{ and } \qquad z_1=-\sqrt{2}+i,
    \qquad z_2 = \sqrt{2} + i. \]
 The free constant is $C =-3/4$. Therefore
\begin{equation} \label{e:Qformula}
    Q(z)=-\frac 14(z+i)^2 (z^2-2i z-3),
\end{equation}
and we recover \eqref{Qr3}.

Once we have $Q(z)$, we may obtain $\mu$ in the following way.
From \eqref{E:Qdefinition} it follows that there is an analytic
branch of $Q^{1/2}(z)$ for $z \in \mathbb C \setminus \gamma$
which behaves as $\tfrac{1}{2} V'(z)$ for large $z$. Choose an orientation
on $\gamma$. The orientation induces a $+$-side and a $-$-side on $\gamma$,
where the $+$-side ($-$-side) is on the left (right) as one traverses the contour according to
its orientation.

\begin{lemma}\label{lemma:probmeasure}
Given the critical trajectory $\gamma$ and the polynomial $Q(z)$, then
\begin{equation} \label{E:muQ}
    \frac{1}{\pi i} \, Q_+^{1/2}(s) \, ds  =  d\mu
\end{equation}
is a probability measure on $\gamma$.
\end{lemma}
\begin{proof}
The measure $\mu$ is a priori complex, however, by the
construction of $\gamma$ we have that
\[ \int_{z_1}^z d\mu(s) = \frac{1}{\pi i} \int_{z_1}^z Q_+^{1/2}(s)
d s \in \mathbb R \]
for every $z \in \gamma$, so that $\mu$ is a real measure.

Taking $z = z_2$ we can compute
\[ \int_{z_1}^{z_2} d\mu(s) = \frac{1}{\pi i} \int_{z_1}^{z_2} Q_+^{1/2}(s) d s \]
by contour integration. Indeed, we have
\[ \frac{1}{\pi i} \int_{z_1}^{z_2} Q_+^{1/2}(s) d s = \frac{1}{2\pi i} \int_C Q^{1/2}(s) d s \]
where $C$ is a closed contour in $\mathbb C \setminus \gamma$ that encircles $\gamma$ once
in the clockwise direction. Moving the contour to infinity, and using
the behavior of $Q^{1/2}$ at infinity, see \eqref{E:Qsqrt}, we find
\begin{equation} \label{E:mumass}
     \mu(\gamma) = \int_{z_1}^{z_2} d\mu(s) = 1.
     \end{equation}

Then if $t \in [0,1] \mapsto z=z(t)$ is a smooth parametrization of $\gamma$ with $z(0)=z_1$
and $z(1) = z_2$ we have that
\[ t \in [0,1] \mapsto \frac{1}{\pi i} \int_{z_1}^{z(t)} Q^{1/2}_+(s) d s \]
is real valued, with values $0$ for $t=0$ and $1$ for $t=1$.
The derivative $z'(t) Q^{1/2}(z(t))$ is non-zero for $0 < t < 1$. Therefore
\[ \frac{1}{\pi i} \int_{z_1}^{z(t)} Q^{1/2}_+(s) d s  \]
is strictly increasing, and it follows that $\mu$ is a probability measure.
\end{proof}

\begin{lemma}\label{lemma:eqmeasure}
Let $\Gamma=\gamma_1\cup\gamma\cup\gamma_2$ be defined as in \eqref{E:Gamma}, \eqref{E:isreal1} and \eqref{E:isreal2}, then the measure $\mu$ defined by \eqref{E:muQ} is the equilibrium measure on $\Gamma$ in the external field $\Re V$.
\end{lemma}
\begin{proof}
From another residue calculation, similar to the one leading to \eqref{E:mumass}
and based on \eqref{E:Qsqrt}, it follows that
\begin{equation} \label{E:muresidue}
    \int_{\gamma} \frac{1}{z-s} d\mu(s) = \frac{1}{2} V'(z) - Q^{1/2}(z),
    \qquad z \in \mathbb C \setminus \gamma.
    \end{equation}

Then we have
\[ g(z) = \int \log(z-s) d\mu(s) \]
is such that \eqref{diffg} holds, which after integration
leads to \eqref{E:gSproperty} and to the first line of
\eqref{E:gequilibrium}.

We extend $\gamma$ to an unbounded contour $\Gamma = \gamma \cup \gamma_1 \cup \gamma_2$
as in section 2.2. The unbounded pieces $\gamma_1$ and $\gamma_2$ are such that
\eqref{E:isreal1} and \eqref{E:isreal2} hold.
This leads to the second line of \eqref{E:gequilibrium}.
For example, if $z \in \gamma_2$, then by \eqref{E:isreal2} and \eqref{E:muresidue}
\begin{align*}
     0 <  2 \int_{z_2}^z Q^{1/2}(s) d s & =
    \int_{z_2}^z (V'(s) - 2 g'(s)) d s \\
    & =
    (V(z) - 2 g(z)) - (V(z_2) - 2g(z_2)) \\
    & = V(z) - 2 g(z) - (\ell + i \tilde{\ell}),
    \end{align*}
which by taking the real part indeed leads to the inequality
in \eqref{E:gequilibrium}.

Because of \eqref{E:gequilibrium} we have that $\mu$ is the
equilibrium measure on $\Gamma$ in the external field $\Re V$.
\end{proof}

Finally, because of \eqref{E:gSproperty} we conclude that the contour $\Gamma$ has the $S$-property,
and this completes the proof of Theorem \ref{theo:Sproperty}.

\section{Proof of Theorem \ref{th:strongasymptotics}} \label{s:proof4}

\subsection{Riemann--Hilbert problem}
The orthogonal polynomial $P_n(z)$ characterized by \eqref{E:varying_weight} appears as the
$(1,1)$ entry of the solution $Y(z)$ of a $2\times 2$ matrix-valued Riemann--Hilbert problem,
see \cite{fokas:1992:isomonodromy}.

From this Riemann--Hilbert problem, the Deift-Zhou steepest descent method performs several explicit
and invertible transformations that allow us to obtain asymptotic results for
the entries of the matrix $Y$, and in particular for $P_n(z)$, as $n\to\infty$ uniformly in
different regions of $\mathbb{C}$, see \cite{DKMVZ:1999:varying}.
In the present case the analysis is quite standard, except for the fact that we are working
on a complex curve $\Gamma$ instead of on a part of the real line. For this reason,
we give a brief sketch of the method and refer the reader to \cite{DKMVZ:1999:varying}, \cite{Deift:2000:RH}
and \cite{KML:2001:Laguerre} for the general theory involving orthogonality with respect
to exponential weights and also for more details  on a similar problem.

We are interested in a matrix-valued function $Y : \mathbb C \setminus \Gamma \to \mathbb{C}^{2\times 2}$ such that
\begin{itemize}
\item $Y(z)$ is analytic for $z \in \mathbb{C}\setminus \Gamma$.
\item $Y_+(z)=Y_-(z) \begin{pmatrix} 1 & e^{-nV(z)} \\ 0 & 1 \end{pmatrix}$, for $z \in \Gamma$,
\item $Y(z)=\left(I+\mathcal{O}\left(\frac 1z\right)\right)
    \begin{pmatrix} z^n & 0 \\ 0 & z^{-n} \end{pmatrix}$, as $z\to\infty$.
\end{itemize}

As before,  $\Gamma$ is the contour $\Gamma = \gamma_1 \cup \gamma \cup \gamma_2$
consisting of the critical trajectory $\gamma$ and its analytic extensions
$\gamma_1$ and $\gamma_2$. See Figure~\ref{Figure3}.

This Riemann--Hilbert problem has a unique solution if and only if
the monic polynomial $P_n(z)$, orthogonal with respect to the
weight function $w(z)$, exists uniquely. If additionally
$P_{n-1}(z)$ exists, then the solution of the Riemann--Hilbert
problem is given by:
\begin{equation*}
Y(z)= \begin{pmatrix} P_n(z)  & (\mathcal{C}P_n w)(z) \\
    -2\pi i\gamma_{n-1}P_{n-1}(z) & -\gamma_{n-1}(\mathcal{C}P_{n-1}w)(z)
    \end{pmatrix},
\end{equation*}
where
\begin{equation}\label{Cauchy}
    (\mathcal{C}f)(z)=\frac{1}{2\pi i}
    \int_{\Gamma}\frac{f(s)}{s-z}\ d s
\end{equation}
is the Cauchy transform on $\Gamma$, and the coefficient
$\gamma_{n-1}$ is defined as
$$
\gamma_{n-1}=\left[\int_{\Gamma}P^2_{n-1}(s)w(s)d s\right]^{-1}.
$$

\subsection{First transformation}
The first transformation $Y\mapsto T$ is a normalization at $\infty$. We
use the functions $\phi_1$ and $\phi_2$ as in \eqref{E:isreal1} and \eqref{E:isreal2}
which are analytic in $\mathbb{C}\setminus(\gamma_1\cup\gamma)$ and $\mathbb{C}\setminus(\gamma_2\cup\gamma)$ respectively and satisfy $\phi_2(z)-\phi_1(z) = \pm \pi i$ for $z\in\mathbb{C}\setminus\Gamma$. We set
\begin{equation}\label{YT}
T(z)=\begin{pmatrix} e^{nl} &0\\0 & e^{-nl}\end{pmatrix}
Y(z)\begin{pmatrix} e^{n[\phi_2(z)-\tfrac 12 V(z)]} & 0\\ 0 &
e^{-n[\phi_2(z)-\tfrac 12 V(z)]}\end{pmatrix}.
\end{equation}

Now, using (\ref{E:Qdefinition}), we obtain by direct integration from \eqref{E:isreal2} that
\begin{equation}
    \phi_2(z)=\tfrac 12 V(z)-\log(z)-l+\mathcal{O}\left(\frac 1z\right), \qquad z\to\infty,
\end{equation}
for some constant of integration $l$. It follows that
$$
e^{n[\phi_2(z)-\tfrac 12 V(z)]}=z^{-n}e^{-nl}\left(1+\mathcal{O}\left(\frac 1z\right)\right), \qquad z\to\infty.
$$

\begin{figure}[t]
\begin{pspicture}(0,0)(10,3.5)
\pscurve{->,arrowsize=0.25}(3,1.5)(4.9,0.8)(6,0.65)
\pscurve(5.9,0.65)(7.1,0.8)(9,1.5)
\psline(1,2.5)(3,1.5)
\psline(9,1.5)(11,2.5)
\psdot(3,1.5)
\psdot(9,1.5)
\put(2.9,1.75){$z_1$}
\put(8.8,1.75){$z_2$}
\put(5.5,0.3){$\gamma$}
\put(1.8,2.3){$\gamma_1$}
\put(9.8,2.3){$\gamma_2$}
\put(9.5,1){$\begin{pmatrix} 1 & e^{-2n\phi_2} \\ 0 & 1 \end{pmatrix}$}
\put(4.5,1.6){$\begin{pmatrix} e^{2n\phi_{2+}} & 1 \\ 0 & e^{2n\phi_{2-}} \end{pmatrix}$}
\put(0.25,1){$\begin{pmatrix} 1 & e^{-2n\phi_1} \\ 0 & 1 \end{pmatrix}$}
\end{pspicture}
\caption{The jump matrices for the Riemann-Hilbert problem for $T$.  \label{Figure6}}
\end{figure}
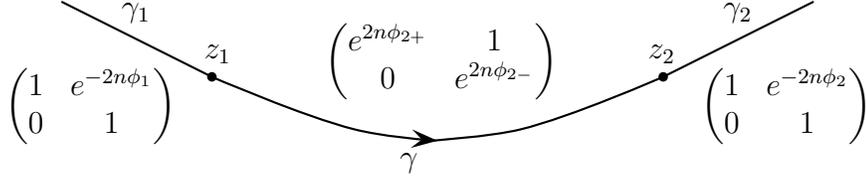

Hence  $T$ satisfies the following Riemann--Hilbert problem:
\begin{itemize}
\item $T(z)$ is analytic for $z$ in $\mathbb{C}\setminus \Gamma$;
\item $T$ has the jumps indicated in Figure \ref{Figure6};
\item $T(z) = I + \mathcal{O}\left(\frac 1z\right)$ as $z\to\infty$.
\end{itemize}

\subsection{Second transformation}

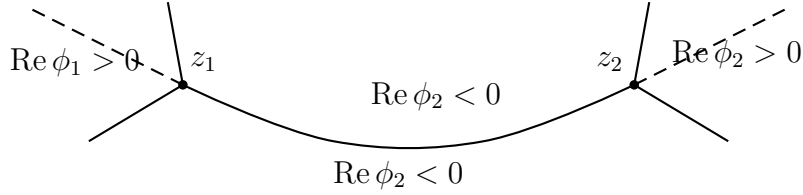
\begin{figure}[ht]
\begin{pspicture}(0,0)(10,3.5)
\pscurve(3,1.5)(5,0.75)(7,0.75)(9,1.5)
\psline[linestyle=dashed](1,2.5)(3,1.5)
\psline[linestyle=dashed](9,1.5)(11,2.5)
\psdot(3,1.5)
\psdot(9,1.5)
\put(3.1,1.75){$z_1$}
\put(8.5,1.75){$z_2$}
\psline(3,1.5)(1.75,0.75)
\psline(3,1.5)(2.8,2.6)
\psline(9,1.5)(10.25,0.75)
\psline(9,1.5)(9.2,2.6)
\put(0.7,1.7){$\Re \phi_1 >0$}
\put(5.5,1.25){$\Re \phi_2<0$}
\put(5,0.2){$\Re \phi_2<0$}
\put(9.5,1.8){$\Re \phi_2>0$}
\end{pspicture}
\caption{The sign of $\Re \phi_1 = \Re \phi_2$ in various parts of the complex plane.
The solid curves are where $\Re \phi_2 = 0$.
The curves $\gamma_1$ and $\gamma_2$ are shown with dashed lines.
We have that $\phi_1$ is real and positive on $\gamma_1$
and $\phi_2$ is real and positive on $\gamma_2$. \label{Figure7}}
\end{figure}

The second transformation of the Riemann-Hilbert problem is
the so-called opening of lenses. From the Cauchy--Riemann equations,
it is possible to show that the sign pattern for  $\Re \phi_2$ is as shown in Figure~\ref{Figure7}.
Since $\phi_1 = \phi_2 \pm \pi i$, the sign pattern for $\Re \phi_1$
is exactly the same.

In the second transformation we open a lens-shaped region around $\gamma$ as in Fig.~\ref{Figure8}, so that the lens is contained in the region where $\Re \phi_2 < 0$:
\begin{figure}[t]
\begin{pspicture}(0,-1)(10,4.5)
\pscurve{->,arrowsize=0.25}(3,1.5)(4.75,0.8)(6,0.65)
\pscurve(5.9,0.65)(7.25,0.8)(9,1.5)
\psline{->,arrowsize=0.25}(0.5,2.75)(1.5,2.25)
\psline(1.5,2.25)(3,1.5)
\psline{->,arrowsize=0.25}(9,1.5)(10.5,2.25)
\psline(10.5,2.25)(11,2.5)
\psdot(3,1.5)
\psdot(9,1.5)
\pscurve{->,arrowsize=0.25}(3,1.5)(3.5,1.75)(4.6,2.15)(6,2.35)
\pscurve(5.9,2.35)(7.4,2.15)(8.5,1.75)(9,1.5)
\pscurve{->,arrowsize=0.25}(3,1.5)(3.5,0.9)(4.5,0.2)(6,-0.15)
\pscurve(5.9,-0.15)(7.5,0.2)(8.5,0.9)(9,1.5)
\put(3,1.75){$z_1$}
\put(8.75,1.75){$z_2$}
\put(9.5,1){$\begin{pmatrix} 1 & e^{-2n\phi_2} \\ 0 & 1 \end{pmatrix}$}
\put(1,3){$\begin{pmatrix} 1 & e^{-2n\phi_1} \\ 0 &  1 \end{pmatrix}$}
\put(5.5,3){$\begin{pmatrix} 1 & 0 \\ e^{2n\phi_2} & 1 \end{pmatrix}$}
\put(5.5,1.25){$\begin{pmatrix} 0 & 1 \\ -1 & 0 \end{pmatrix}$}
\put(1.75,-0.25){$\begin{pmatrix} 1 & 0 \\ e^{2n\phi_2} & 1  \end{pmatrix}$}
\put(1,1.85){$\gamma_1$}
\put(10,2.5){$\gamma_2$}
\put(4.5,1.1){$\gamma$}
\end{pspicture}
\caption{The contour $\Gamma_S$
and the jump matrices on $\Gamma_S$ in  the Riemann-Hilbert problem for $S$.
\label{Figure8}}
\end{figure}
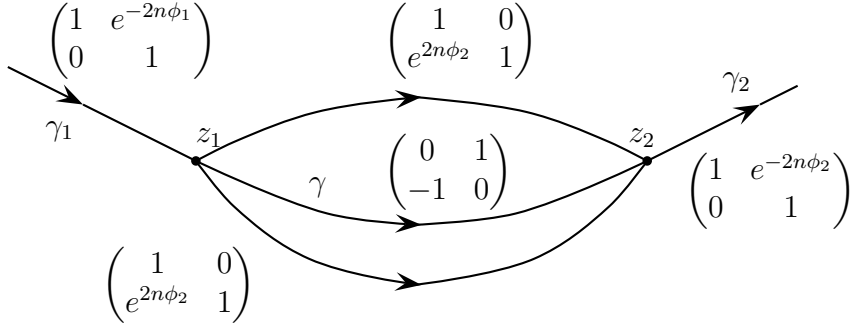

We define
\begin{equation}
\label{matrixS}
S =  \begin{cases}
   \, T \begin{pmatrix} 1 & 0 \\ -e^{2n\phi_2} & 1 \end{pmatrix} & \text{in the upper part of the lens}, \\
   \,T \begin{pmatrix}  1 & 0 \\ e^{2n\phi_2} & 1 \end{pmatrix} & \text{in the lower part of the lens}, \\
   \, T & \text{elsewhere}.
    \end{cases}
\end{equation}
Then $S$ satisfies the following Riemann--Hilbert problem:
\begin{itemize}
\item $S(z)$ is analytic for $z \in \mathbb{C}\setminus \Gamma_S$,
    where $\Gamma_S$ consists of $\Gamma$ plus the lips of the lens;
\item $S$ has the jumps indicated in Figure \ref{Figure8};
\item $S(z) = I +\mathcal{O}\left(\frac 1z\right)$ as $z\to\infty$.
\end{itemize}

\subsection{Construction of parametrices}
\subsubsection{Global parametrix}
Now we seek an approximation to $S$ that is valid for large $n$.
The approximation will consist of two parts, a global parametrix
$N$ away from the endpoints $z_1$ and $z_2$ and local parametrices $P$ at $z_1$ and $z_2$.

The global parametrix $N$ satisfies a Riemann--Hilbert problem
with the same constant jump on $\gamma$. Then $R=SN^{-1}$ will be
analytic across $\gamma$. We define $N$ as
\begin{equation} \label{E:defN}
    N(z) =
    \begin{pmatrix} \frac 12(\beta(z)+\beta(z)^{-1}) & \frac {1}{2i}(\beta(z)-\beta(z)^{-1}) \\[5pt]
     -\frac {1}{2i}(\beta(z)-\beta(z)^{-1}) & \frac 12(\beta(z)+\beta(z)^{-1}) \end{pmatrix},
\end{equation}
where $\beta$ is given by \eqref{E:beta}, see \cite{Deift:2000:RH}, \cite{DKMVZ:1999:varying} and \cite{Arno:2003:RH}.

Then $N$ satisfies the following Riemann--Hilbert problem:
\begin{itemize}
\item $N(z)$ is analytic for $z \in \mathbb{C}\setminus \gamma$;
\item $N_+ = N_- \begin{pmatrix} 0 & 1 \\ -1 & 0 \end{pmatrix}$ on $\gamma$;
\item $N(z)=I+\mathcal{O}\left(\frac 1z\right)$ as $z\to\infty$;
\item $N(z)=\mathcal{O}(|z-z_j|^{-1/4})$ as $z\to z_j$ for $j=1,2$.
\end{itemize}
It is clear that $N$ cannot be a good approximation to $S$ near the endpoints of $\gamma$,
since it blows up at $z=z_1$ and $z=z_2$, while $S$ remains bounded there.
For this reason we need a different local approximation near the endpoints.

\subsection{Local parametrix}
The local parametrix $P$ is constructed in neighbourhoods of the endpoint $z=z_j$, $j=1,2$, say
$$
U_{\delta}(z_j)=\{z\in\mathbb{C} \mid |z-z_j|<\delta\}, \qquad j =1,2,
$$
with some small but fixed $\delta > 0$. We describe here the
construction of $P$ in $U_{\delta}(z_2)$, the construction in $U_{\delta}(z_1)$
being similar.

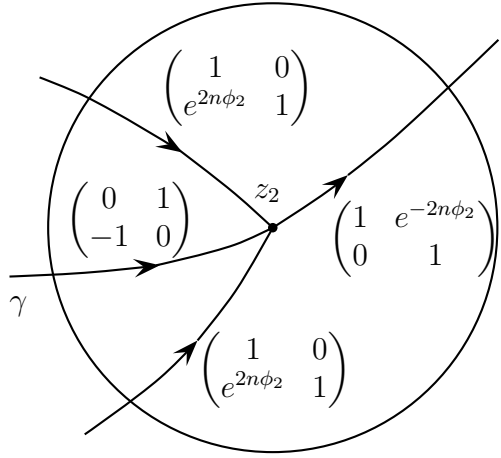
\begin{figure}[t]
\begin{pspicture}(0,0)(10,6.5)
\psdot(7,3)
\pscircle(7,3){3}
\pscurve{->,arrowsize=0.25}(3.5,2.35)(4.5,2.4)(5.5,2.5)
\pscurve(5.45,2.5)(6.5,2.775)(7,3)
\pscurve{->,arrowsize=0.25}(3.9,5)(4.5,4.75)(5.8,4)
\pscurve(5.7,4.075)(6.5,3.45)(7,3)
\pscurve{->,arrowsize=0.25}(4.5,0.25)(5.5,1)(6,1.5)
\pscurve(6,1.5)(6.5,2.15)(7,3)
\pscurve{->,arrowsize=0.25}(7,3)(7.75,3.5)(8,3.7)
\pscurve(8,3.7)(8.5,4.1)(9,4.55)(10,5.5)
\put(6.75,3.4){$z_2$}
\put(3.5,1.9){$\gamma$}
\put(4.25,3){$\begin{pmatrix} 0 & 1  \\ -1 & 0 \end{pmatrix}$}
\put(7.75,2.75){$\begin{pmatrix} 1 & e^{-2n\phi_2} \\ 0 & 1 \end{pmatrix}$}
\put(5.5,4.75){$\begin{pmatrix} 1 & 0 \\ e^{2n\phi_2} & 1\end{pmatrix}$}
\put(6,1){$\begin{pmatrix} 1 & 0 \\ e^{2n\phi_2} & 1 \end{pmatrix}$}
\end{pspicture}
\caption{The jump matrices in the Riemann-Hilbert problem for $P$
defined in the neighbourhood $U_{\delta}(z_2)$ of $z_2$. \label{Figure9}}
\end{figure}

The local parametrix $P$
should  satisfy the following Riemann--Hilbert problem:
\begin{itemize}
\item $P(z)$ is analytic for $z \in U_{\delta}(z_2)\setminus \Gamma_S$
    with a continuous extension to $\overline{U_{\delta}(z_2)} \setminus \Gamma_S$;
\item $P$ has the jumps on $\Gamma_S \cap U_{\delta}(z_2)$ as shown in Fig.~\ref{Figure9} (these are the same jump
    matrices as in the RH problem for $S$);
\item $P(z) = \left(I+\mathcal{O}\left(\frac 1n\right)\right)N(z)$ as $n\to\infty$,
    uniformly for $z \in \partial U_{\delta}(z_2)$;
\item $P(z)$ remains bounded as $z\to z_2$.
\end{itemize}

The construction of $P$ is given in terms of the Airy function $\Ai$
and its derivative.  We put
\begin{equation}  \label{E:defP}
   P(z) = E_n(z)A(n^{2/3} f(z)) \begin{pmatrix} e^{n\phi_2(z)} & 0 \\ 0 & e^{-n\phi_2(z)} \end{pmatrix},
\end{equation}
where $A(\zeta)$, $f(z)$, and $E_n(z)$ are described below.

\begin{figure}[t]
\begin{pspicture}(0,0)(10,7)
\psline{->,arrowsize=0.25}(3.25,6)(4.125,4.5)
\psline(4.125,4.5)(5,3)
\psline{->,arrowsize=0.25}(3.15,0)(4.125,1.5)
\psline(4.125,1.5)(5,3)
\psline{->,arrowsize=0.25}(5,3)(7,3)
\psline(7,3)(9,3)
\psline{->,arrowsize=0.25}(1,3)(3,3)
\psline(3,3)(5,3)
\psarc(5,3){1}{0}{120}
\put(7,3.5){$\begin{pmatrix} 1 & 1 \\ 0 & 1 \end{pmatrix}$}
\put(4.25,5){$\begin{pmatrix} 1 & 0 \\ 1 & 1 \end{pmatrix}$}
\put(4.25,1){$\begin{pmatrix} 1 & 0 \\ 1 & 1 \end{pmatrix}$}
\put(2,3.5){$\begin{pmatrix} 0 & 1 \\ -1 & 0 \end{pmatrix}$}
\put(5.75,3.8){$2\pi/3$}
\end{pspicture}
\caption{Contours and jump matrices in the Riemann-Hilbert problem
for Airy functions. \label{Figure10}}
\end{figure}
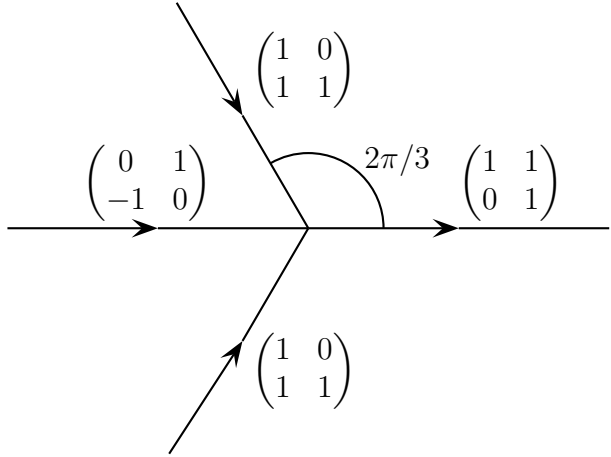

\paragraph{Airy parametrix $A(\zeta)$}
The matrix-valued function $A(\zeta)$ is the solution
of the Airy Riemann-Hilbert problem, which is posed on
four infinite rays in an auxiliary $\zeta$-plane as follows:
\begin{itemize}
\item $A(\zeta)$ is analytic for $\zeta \in \mathbb C$, $\arg \zeta \not\in \{0, 2\pi/3, -2\pi/3, \pi\}$;
\item $A$ has the jumps on the four rays as shown in Fig.~\ref{Figure10};
\item As $\zeta \to \infty$, we have
\begin{equation}\label{asympA}
A(\zeta)=  \begin{pmatrix} \zeta^{-1/4} & 0 \\ 0 & \zeta^{1/4} \end{pmatrix}
\frac{1}{\sqrt{2}} \begin{pmatrix} 1 & i \\ i & 1 \end{pmatrix}
     \left(I+\mathcal{O}\left(\frac{1}{\zeta^{3/2}}\right)\right)
     \begin{pmatrix} e^{-\frac 23 \zeta^{3/2}} & 0  \\ 0 & e^{\frac 23 \zeta^{3/2}} \end{pmatrix}
\end{equation}
\item $A(\zeta)$ remains bounded as $\zeta \to 0$.
\end{itemize}

The solution of this Riemann--Hilbert problem is given by
the Airy function $\Ai(\zeta)$ and rotated versions of it, see \cite{DKMVZ:1999:varying}
and \cite[Sec.~10.4]{Abramowitz:1964:HMF}. Let
\begin{equation*}
y_0(\zeta)=\Ai(\zeta), \quad y_1(\zeta)=\omega\Ai(\omega\zeta), \quad y_2(\zeta)=\omega^2\Ai(\omega^2\zeta),
\end{equation*}
where $\omega=e^{\tfrac{2\pi i}{3}}$. These are three solutions of the Airy differential
equation $y'' = \zeta y$ satisfying the connection formula
$y_0(\zeta)+y_1(\zeta)+y_2(\zeta)=0$.
For instance, in the sector $0<\arg\zeta<2\pi/3$ we have
\begin{equation} \label{Ainfirstsector}
A(\zeta)=\sqrt{2\pi} \begin{pmatrix} y_0(\zeta) & -y_2(\zeta) \\
    -i y_0'(\zeta) & i y_2'(\zeta) \end{pmatrix}.
\end{equation}
The
solution in the other sectors is obtained from this by applying
the appropriate jump matrices.

\paragraph{Conformal map $f(z)$}
The map $f(z)$ is defined by
\begin{equation}
    f(z)= \left[\frac {3}{2}\phi_2(z)\right]^{2/3},
\end{equation}
which is a conformal map in a neighourhood of $z=z_2$. It is assumed that $\delta > 0$
is sufficiently small so that $f$ is indeed a conformal map on $U_{\delta}(z_2)$, and also that the lens around $\gamma$ is opened in such a way that the
lips of the lens inside $U_{\delta}(z_2)$ are mapped by $\zeta = f(z)$ to
the rays $\arg \zeta = \pm 2\pi/3$. This can be done without any loss of generality.

\paragraph{Analytic prefactor $E_n(z)$}
The prefactor $E_n(z)$ in \eqref{E:defP} is defined by
\begin{align}
    E_n(z) & = \nonumber
    N(z)\frac{1}{\sqrt{2}} \begin{pmatrix} 1 & - i \\ -i & 1 \end{pmatrix}
    \begin{pmatrix} n^{1/6} f(z)^{1/4} & 0 \\ 0 & n^{-1/6} f(z)^{-1/4} \end{pmatrix} \\
    & = \frac{1}{\sqrt{2}} \begin{pmatrix} 1 & - i \\ -i & 1 \end{pmatrix}
        \begin{pmatrix} n^{1/6} f(z)^{1/4} \beta^{-1}(z) & 0 \\
        0 & n^{-1/6} f(z)^{-1/4} \beta(z) \end{pmatrix}, \label{matrixEn}
\end{align}
which is analytic in $U_{\delta}(z_2)$.
It is chosen so that the matching condition $P(z) = (I + \mathcal{O}(1/n)) N(z)$
for $z \in \partial U_{\delta}(z_2)$ is satisfied.

Then with these definitions it can be shown that $P$ defined by \eqref{E:defP}
indeed satisfies the Riemann-Hilbert problem for $P$.

\subsection{Third transformation}

In the third and final transformation we use the global parametrix $N(z)$
and the local parametrices $P(z)$ to define
\begin{align} \label{matrixR}
R(z) = \begin{cases} S(z)N(z)^{-1}, & \quad z\in\mathbb{C} \setminus (\Gamma_S \cup
    \overline{U_{\delta}(z_1)} \cup \overline{U_{\delta}(z_2)}), \\[5pt]
    S(z)P(z)^{-1},  & \quad z\in(\overline{U_{\delta}(z_1)}
    \cup \overline{U_{\delta}(z_2)}) \setminus \Gamma_S.
    \end{cases}
\end{align}

\begin{figure}[t]
\begin{pspicture}(0,0)(15,7)
\psdot(3,3.5)
\psdot(9,3.5)
\pscurve{->,arrowsize=0.25}(3.8,4.1)(4.5,4.3)(5,4.4)(6,4.5)
\pscurve(6,4.5)(7,4.4)(7.5,4.3)(8.2,4.1)
\pscurve{->,arrowsize=0.25}(3.65,2.75)(4.5,2.3)(5,2.15)(6,2)
\pscurve(6,2)(7,2.15)(7.5,2.3)(8.35,2.75)
\put(2.8,3.75){$z_1$}
\put(8.8,3.75){$z_2$}
\psarc(9,3.5){1}{-90}{90}
\psarc{<-,arrowsize=0.25}(9,3.5){1}{90}{-90}
\psarc(3,3.5){1}{-90}{90}
\psarc{<-,arrowsize=0.25}(3,3.5){1}{90}{-90}
\pscurve{->,arrowsize=0.25}(9.95,3.75)(10.5,4)(11,4.25)
\pscurve(11,4.25)(11.5,4.6)(12,5)
\pscurve(2.05,3.75)(1.5,4)(1,4.25)
\pscurve{<-,arrowsize=0.25}(1,4.25)(0.5,4.6)(0,5)
\put(4.5,5){$N\begin{pmatrix} 1 & 0 \\ e^{2n\phi_2} & 1 \end{pmatrix} N^{-1}$}
\put(4.5,1){$N\begin{pmatrix} 1 & 0 \\ e^{2n\phi_2} & 1 \end{pmatrix} N^{-1}$}
\put(10,2.75){$PN^{-1}$}
\put(1,2.75){$PN^{-1}$}
\put(8.5,5.5){$N\begin{pmatrix} 1 & e^{-2n\phi_2} \\ 0 & 1 \end{pmatrix} N^{-1}$}
\put(0,5.5){$N\begin{pmatrix} 1 & e^{-2n\phi_1} \\ 0 & 1 \end{pmatrix} N^{-1}$}
\end{pspicture}
\caption{The contour $\Gamma_R$
and the jump matrices on $\Gamma_R$ in  the Riemann-Hilbert problem for $R$.
\label{Figure11}}
\end{figure}
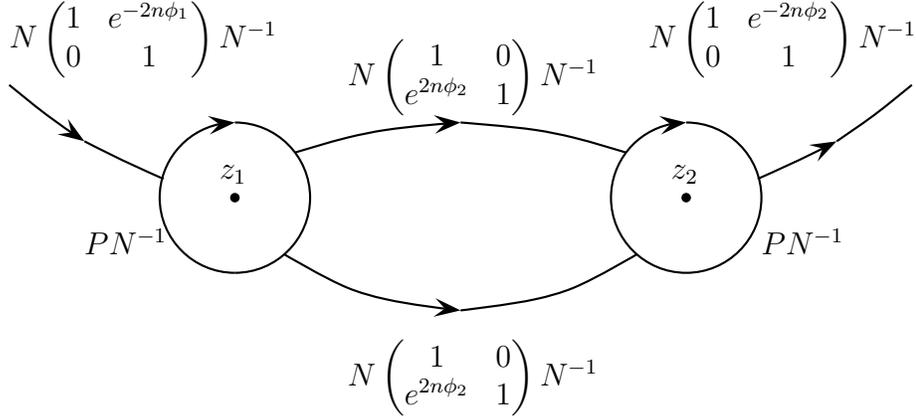

Then $R$ has an analytic continuation across $\gamma$ and across the
parts of $\Sigma_S$ that are inside the disks $U_{\delta}(z_1)$
and $U_{\delta}(z_2)$. It satisfies the following Riemann--Hilbert problem:

\begin{itemize}
\item $R(z)$ is analytic for $z \in \mathbb{C}\setminus \Gamma_R$, where $\Gamma_R$ is
    the contour shown in Fig.~\ref{Figure11};
\item $R$ has jumps on each part of $\Gamma_R$ with jump matrices
    as indicated in Fig.~\ref{Figure11};
\item $R(z)=I+\mathcal{O}\left(\frac 1z\right)$ as $z\to\infty$.
\end{itemize}

The jump matrices in the Riemann-Hilbert problem for $R$ tend
to the identity matrix as $n \to \infty$.
Indeed, since $P(z) = (I+\mathcal{O}(1/n))N(z)$ as $n\to\infty$, uniformly
for $z \in \partial U_{\delta}(z_1)\cup\partial U_{\delta}(z_2)$, we have that
\[
    R_+(z)=R_-(z)\left(I+\mathcal{O}\left(\frac 1n\right)\right),
    \qquad z\in \partial U_{\delta}(z_1)\cup\partial U_{\delta}(z_2)
\]
as $n \to \infty$. On the remaining parts of $\Gamma_R$ we even have
for some positive constant $c > 0$,
\[
    R_+(z) = R_-(z)\left(I+\mathcal{O}(e^{-cn})\right),
     \qquad z\in \Gamma_R\setminus(\partial U_{\delta}(z_1)\cup\partial U_{\delta}(z_2)),
\]
as $n \to \infty$. Thus the jumps on $R$ tend to the identity matrix uniformly,
and in fact also in $L^2(\Gamma_R)$.

It then follows from the general theory, see \cite{Deift:2000:RH}
and \cite{Arno:2003:RH}, that
the solution to the Riemann-Hilbert problem for $R$ exists for all large
enough $n$ with
\begin{equation} \label{E:Rasymptotics}
    R(z)=I+\mathcal{O}\left(\frac{1}{n}\right), \qquad \text{as } n\to\infty,
\end{equation}
uniformly for $z\in\mathbb{C}\setminus\Gamma_R$.

\subsection{Proof of Theorem \ref{th:strongasymptotics}}
Once we arrive at this result for $R$, it is possible to reverse
all the transformations $Y \mapsto T \mapsto S \mapsto R$, since they are all explicit and
invertible. The first thing that follows is that the original Riemann-Hilbert problem for $Y$
has a unique solution for large enough $n$.
Since
\[ P_n(z) = Y_{11}(z), \]
this proves that the orthogonal polynomials $P_n$
indeed exist for every large enough $n$.

The asymptotic formula \eqref{E:Rasymptotics} for $R$  further yields the
first term in an asymptotic expansion of $Y$ as $n \to \infty$.
Following the effect of the inverse transformations $R \mapsto S \mapsto T \mapsto Y$
on the asymptotic formula  \eqref{E:Rasymptotics} for $R$,
we obtain the asymptotics of $Y$ and therefore of $P_n$ in the various
regions of the complex plane.
This will give the different parts of Theorem \ref{th:strongasymptotics}.

\subsubsection{Proof of part (a)}
Let $z \in \mathbb C \setminus \gamma$. We then may and do assume
that the lens around $\gamma$ and the neighborhoods $U_{\delta}(z_1)$
and $U_{\delta}(z_2)$ are chosen so that $z$ is in the outside region.

From \eqref{matrixR} we have that $S(z) = R(z) N(z)$.
Also, because we are outside the lens, we have $S(z) = T(z)$ from (\ref{matrixS})
and $T(z)$ in terms of $Y(z)$ follows from (\ref{YT}). Combining all this we find
\[ Y(z) = \begin{pmatrix} e^{-nl} &  0 \\ 0 & e^{nl} \end{pmatrix}
    \left(I + \mathcal O\left(\frac{1}{n}\right)\right) N(z)
    \begin{pmatrix} e^{-n\left[\phi_2(z)-\frac{1}{2} V(z)\right]} & 0 \\
        0 & e^{n\left[\phi_2(z)-\frac{1}{2} V(z)\right]} \end{pmatrix}.
\]
Then, for the $(1,1)$-entry the first part of the theorem follows
in a straightforward way, since
\[ g(z) = \frac{1}{2} V(z)  - \phi_2(z) - l. \]

\subsubsection{Proof of part (b)}

For $z$ inside the lens, but outside of the two disks, we
have $S(z) = R(z) N(z)$ as before. From \eqref{matrixS} we
then get
\begin{equation*}
    T(z)=S(z) \begin{pmatrix} 1 & 0 \\ \pm e^{2n\phi_2} & 1 \end{pmatrix} =
        R(z)N(z) \begin{pmatrix} 1 & 0 \\ \pm e^{2n\phi_2} & 1 \end{pmatrix}
\end{equation*}
where the $+$ sign ($-$ sign) is taken in the upper (lower) part of the lens.
Using \eqref{YT} and \eqref{E:Rasymptotics} we then find
\begin{equation*}
    Y(z)= \begin{pmatrix} e^{-nl} & 0 \\ 0 & e^{nl} \end{pmatrix} \left(I + \mathcal O\left(\frac{1}{n}\right)\right) N(z)
    \begin{pmatrix} e^{-n\left[\phi_2(z)-\frac{1}{2} V(z)\right]} & 0 \\
        \pm e^{n\left[\phi_2(z)+\frac{1}{2} V(z) \right]} & e^{n\left[\phi_2(z)-\frac{1}{2} V(z)\right]} \end{pmatrix}.
\end{equation*}
Then for the $(1,1)$-entry we obtain from this
\begin{align*}
    P_n(z) & = Y_{11}(z)  = e^{n \left[\frac{V(z)}{2} - l\right]}
        \begin{pmatrix} 1 + \mathcal O\left(\frac{1}{n}\right), & \mathcal O\left(\frac{1}{n}\right) \end{pmatrix}
            N(z) \begin{pmatrix} e^{-n \phi_2(z)} \\ \pm e^{n \phi_2(z)} \end{pmatrix} \\
            & = e^{n\left[\frac {V(z)}{2}-l\right]}
            \left( e^{-n \phi_2(z)} N_{11}(z) \pm e^{n\phi_2(z)} N_{12}(z) + \mathcal O\left(\frac{1}{n}\right)
            \right)
\end{align*}
as $n\to\infty$. This proves part (b) of the theorem.

\subsection{Proof of part (c)}

In the neighbourhoods $U_{\delta}(z_1)$ and $U_{\delta}(z_2)$ of
the endpoints $z_1$ and $z_2$ we use the local parametrix $P(z)$
to obtain  an approximation for $P_n(z)$ in terms of Airy
functions. Indeed, by \eqref{matrixR} and \eqref{E:Rasymptotics},
\begin{equation*}
    S(z) = R(z) P(z) = \left( I + \mathcal O\left(\frac{1}{n}\right)\right) P(z)
\end{equation*}
for $z \in U_{\delta}(z_1) \cup U_{\delta}(z_2)$. If we assume that
$z$ is inside the disk $U_{\delta}(z_2)$ but outside the lens around $\gamma$,
then we find by following the transformations \eqref{matrixS} and \eqref{YT} that
\begin{align*}
    Y(z)= \begin{pmatrix} e^{-nl} & 0 \\ 0 & e^{nl} \end{pmatrix}
        \left( I + \mathcal O\left(\frac{1}{n}\right)\right)
    P(z) \begin{pmatrix} e^{-n\left[\phi_2(z)-\frac{V(z)}{2}\right]} & 0 \\
        0  & e^{n\left[\phi_2(z)-\frac{V(z)}{2}\right]} \end{pmatrix}.
\end{align*}

Using \eqref{E:defP} and \eqref{matrixEn}, we obtain from this that
\begin{multline*}
Y(z)= \frac{1}{\sqrt{2}} \begin{pmatrix} e^{-nl} & 0 \\ 0 & e^{nl} \end{pmatrix}
    \left( I + \mathcal O\left(\frac{1}{n}\right)\right)  \begin{pmatrix} 1 & -i \\ -i & 1 \end{pmatrix} \\
    \times  \begin{pmatrix} n^{1/6} f(z)^{1/4} \beta(z)^{-1} & 0 \\
        0  & n^{-1/6} f(z)^{-1/4} \beta(z) \end{pmatrix} \\
    \times A(n^{2/3} f(z))
    \begin{pmatrix} e^{n\frac{V(z)}{2}} & 0  \\ 0 & e^{-n\frac{V(z)}{2}} \end{pmatrix}.
\end{multline*}

To evaluate $A(n^{2/3} f(z))$ we use \eqref{Ainfirstsector} and it follows that
\begin{align*}
     P_{n}(z) & = \begin{pmatrix} 1 & 0 \end{pmatrix}  Y(z) \begin{pmatrix} 1 \\ 0 \end{pmatrix} \\
     &=
    \sqrt{\pi} e^{n \left[\frac{V(z)}{2} - l\right]}
    \begin{pmatrix} 1 + \mathcal O\left(\frac{1}{n} \right), & \mathcal O\left(\frac{1}{n}\right) \end{pmatrix}
    \begin{pmatrix} 1 & -i \\ -i & 1 \end{pmatrix} \\
    & \qquad \times \begin{pmatrix} n^{1/6} f(z)^{1/4} \beta(z)^{-1} & 0  \\
        0 & n^{-1/6} f(z)^{-1/4} \beta(z) \end{pmatrix}
    \begin{pmatrix} \Ai(n^{2/3} f(z)) \\ - i \Ai'(n^{2/3} f(z)) \end{pmatrix}
\end{align*}
as $n\to\infty$.
This proves part (c) of the theorem in case $z \in U_{\delta}(z_2)$ is outside the lens.
A similar calculation leads to the same expression in case $z$ is inside the lens.
This completes the proof of part (c) of Theorem \ref{th:strongasymptotics}.

\section{Concluding remarks}
We have presented a Riemann--Hilbert analysis of a family of
polynomials orthogonal with respect to a varying exponential
weight on certain curves of the complex plane. The problem was
motivated by the fact that the zeros of these polynomials are
complex Gaussian quadrature points for an oscillatory integral on
an interval $[a,b]\subset\mathbb{R}$. The zeros cluster on
analytic arcs in the complex plane, which are given by a critical
trajectory of a suitable quadratic differential.

We have focused on the case where the weight function is
$V(z)= - i z^3/3$, for which we were able to obtain explicit expressions
throughout the Riemann-Hilbert analysis. A similar procedure (with
more complicated computations) can be applied in principle to the more general
case $V(z)= - i z^r/r$ with $r\geq 5$ and odd. The only difficulty
is the determination of a curve with the $S$-property in this
more general case. It would be interesting to know
if we are in the one-cut case for every odd $r$.

It is also worth remarking that the Riemann--Hilbert analysis can
provide  more detailed asymptotic information than the one given
before, following the ideas exposed in  \cite{DKMVZ:1999:varying}, \cite{KMVV:2004:Jacobi}. The
importance of these results from a numerical point of view is
currently under investigation.

\section*{Acknowledgements}
The authors acknowledge useful discussions with
A.~Mart\'{i}nez-Finkelshtein and H.~Stahl. A. Dea\~{n}o
acknowledges financial support from the programme of postdoctoral
grants of the Spanish Ministry of Education and Science and
project MTM2006-09050. D. Huybrechs is a Postdoctoral Fellow of the Research Foundation
Flanders (FWO) and is supported by FWO-Flanders project G061710N.
A.B.J.~Kuijlaars
is supported by K.U. Leuven research grant OT/08/33, FWO-Flanders project
G.0427.09, by  the Belgian Interuniversity Attraction Pole P06/02, by the
European Science Foundation Program MISGAM, and by grant
MTM2008-06689-C02-01 of the Spanish
Ministry of Science and Innovation.

\bibliography{RH}

\begin{thebibliography}{10}
\expandafter\ifx\csname url\endcsname\relax
  \def\url#1{\texttt{#1}}\fi
\expandafter\ifx\csname urlprefix\endcsname\relax\def\urlprefix{URL }\fi

\bibitem{Abramowitz:1964:HMF}
M.~Abramowitz, I.~A. Stegun, Handbook of Mathematical Functions, vol.~55 of
  National Bureau of Standards Applied Mathematics Series, U.S. Government
  Printing Office, Washington, 1964.

\bibitem{Apt:2002}
A.~Aptekarev, Sharp constants for rational approximation of analytic functions,
  Sbornik Math. 193 (2002) 1--72.

\bibitem{Apt:2007:complex}
A.~Aptekarev, R.~Khabibullin, Asymptotic expansions for polynomials orthogonal
  with respect to a complex non-constant weight function, Trans. Moscow Math.
  Soc. 68 (2007) 1--37.

\bibitem{Baik:2001:random}
J.~Baik, P.~Deift, K.~T.-R. McLaughlin, P.~Miller, X.~Zhou, Optimal tail
  estimates for directed last passage site percolation with geometric random
  variables, Adv. Theor. Math. Phys. 5 (2001) 1207--1250.

\bibitem{Bertola:Lp}
F.~Balogh, M.~Bertola, On the norms and roots of orthogonal polynomials in the
  plane and ${L}_p$-optimal polynomials with respect to varying weights,
  arXiv:0910.4223v1.

\bibitem{Bertola:Boutroux}
M.~Bertola, Boutroux curves with external field: equilibrium measures without a
  variational problem, arXiv:0705.3062.

\bibitem{BertolaMo:2009}
M.~Bertola, M.~Y. Mo, Commuting difference operators, spinor bundles and the
  asymptotics of orthogonal polynomials with respect to varying complex
  weights, Adv. Math. 220 (2009) 154--218.

\bibitem{DH:2008:CG}
A.~Dea\~{n}o, D.~Huybrechs, {C}omplex {G}aussian quadrature of oscillatory
  integrals, Numer. Math. 112~(2) (2009) 197--219.

\bibitem{Deift:2000:RH}
P.~Deift, Orthogonal Polynomials and Random Matrices: a {R}iemann--{H}ilbert
  Approach, American Mathematical Society, Providence, RI, 1999.

\bibitem{DKMVZ:1999:varying}
P.~Deift, T.~Kriecherbauer, K.~T.-R. McLaughlin, S.~Venakides, X.~Zhou, Uniform
  asymptotics for polynomials orthogonal with respect to varying exponential
  weights and applications to universality questions in random matrix theory,
  Comm. Pure Appl. Math. 52 (1999) 1335--1425.

\bibitem{deift:1993:steepestdescent}
P.~Deift, X.~Zhou, A steepest descent method for oscillatory
  {R}iemann-{H}ilbert problems. {A}symptotics for the {MKdV} equation, Ann.
  Math. 137 (1993) 295--368.

\bibitem{fokas:1992:isomonodromy}
A.~S. Fokas, A.~R. Its, A.~V. Kitaev, The isomonodromy approach to matrix
  models in 2{D} quantum gravity, Comm. Math. Phys. 147 (1992) 395--430.

\bibitem{GR:1989:eq}
A.~A. Gonchar, E.~A. Rakhmanov, Equilibrium distributions and degree of
  rational approximation of anaytic functions, Math. USSR Sbornik 62 (1989)
  305--348.

\bibitem{huybrechs:2009:hoq}
D.~Huybrechs, S.~Olver, Highly oscillatory quadrature, in: B.~Engquist,
  A.~Fokas, E.~Hairer, A.~Iserles (eds.), Highly Oscillatory Problems,
  Cambridge Univ. Press, Cambridge, 2009, pp. 25--50.

\bibitem{Huybrechs:06:osc1}
D.~Huybrechs, S.~Vandewalle, On the evaluation of highly oscillatory integrals
  by analytic continuation, SIAM J. Numer. Anal. 44~(3) (2006) 1026--1048.

\bibitem{Arno:2003:RH}
A.~B.~J. Kuijlaars, {R}iemann-{H}ilbert analysis for orthogonal polynomials,
  in: E.~Koelink, W.~Van~Assche (eds.), Orthogonal Polynomials and Special
  Functions, vol. 1817 of Lecture Notes in Mathematics, Springer Verlag,
  Berlin, 2003, pp. 167--210.

\bibitem{KMF:2004:Jacobi}
A.~B.~J. Kuijlaars, A.~Mart\'{i}nez-Finkelshtein, Strong asymptotics for
  {J}acobi polynomials with varying nonstandard parameters, J. Anal. Math. 94
  (2004) 195--234.

\bibitem{KML:2001:Laguerre}
A.~B.~J. Kuijlaars, K.~T.-R. McLaughlin, Riemann-{H}ilbert analysis for
  {L}aguerre polynomials with large negative parameter, Comput. Meth. Funct.
  Theory 1 (2001) 205--233.

\bibitem{KML:2004:Laguerre}
A.~B.~J. Kuijlaars, K.~T.-R. McLaughlin, Asymptotic zero behavior of {L}aguerre
  polynomials with negative parameter, Constr. Approx. 20 (2004) 497--523.

\bibitem{KMVV:2004:Jacobi}
A.~B.~J. Kuijlaars, K.~T.-R. McLaughlin, W.~Van~Assche, M.~Vanlessen, The
  {R}iemann-{H}ilbert approach to strong asymptotics for orthogonal polynomials
  on $[-1,1]$, Adv. Math. 188~(2) (2004) 337--398.

\bibitem{MF:1997}
A.~Mart\'{\i}nez-Finkelshtein, {T}rajectories of quadratic differentials and
  approximations of exponents on the semiaxis, in:
  A.~Mart\'{\i}nez-Finkelshtein, F.~Marcell\'an, J.~Moreno (eds.), Complex
  Methods in Approximation Theory, Universidad de Almer\'{\i}a, 1997, pp.
  69--84.

\bibitem{MFMGO:2001:Laguerre}
A.~Mart\'{i}nez-Finkelshtein, P.~Mart\'{i}nez-Gonz\'{a}lez, R.~Orive, On
  asymptotic zero distribution of {L}aguerre and generalized {B}essel
  polynomials with varying parameters, J. Comp. Appl. Math. 133~(1-2) (2001)
  477--487.

\bibitem{MFO:2005:Jacobi}
A.~Mart\'{i}nez-Finkelshtein, R.~Orive, Riemann-{Hilbert} analysis for {Jacobi}
  polynomials orthogonal on a single contour, J. Approx. Theory 134 (2005)
  137--170.

\bibitem{ST:1997:Pot}
E.~Saff, V.~Totik, Logarithmic Potentials with External Fields, Springer
  Verlag, Berlin, 1997.

\bibitem{Stahl:1986}
H.~Stahl, Orthogonal polynomials with complex-valued weight function. {I},
  {II}, Constr. Approx. (1986) 225--240, 241--251.

\bibitem{strebel:1984:quad}
K.~Strebel, Quadratic Differentials, Springer Verlag, Berlin, 1984.

\bibitem{wong:2001:asymptotic}
R.~Wong, Asymptotic Approximation of Integrals, SIAM, Philadelphia, 2001.

\end{thebibliography}
\bibliographystyle{elsart-num-sort}

\end{document}